# Numerical Solution of second order hyperbolic telegraph equation via new Cubic Trigonometric B-Splines Approach

Tahir Nazir[a], Muhammad Abbas[a,*], Ahmad Izani Md. Ismail[b], Ahmad Abd. Majid[b]
[a]Department of Mathematics, University of Sargodha, 40100 Sargodha, Pakistan
[b]School of Mathematical Sciences, Universiti Sains Malaysia, 11800 Penang Malaysia

**Abstract**
This paper presents a new approach and methodology to solve the second order one dimensional hyperbolic telegraph equation with Dirichlet and Neumann boundary conditions using the cubic trigonometric B-spline collocation method. The usual finite difference scheme is used to discretize the time derivative. The cubic trigonometric B-spline basis functions are utilized as an interpolating function in the space dimension, with a $\theta$ weighted scheme. The scheme is shown to be unconditionally stable for a range of $\theta$ values using the von Neumann (Fourier) method. Several test problems are presented to confirm the accuracy of the new scheme and to show the performance of trigonometric basis functions. The proposed scheme is also computationally economical and can be used to solve complex problems. The numerical results are found to be in good agreement with known exact solutions and also with earlier studies.

**Keywords:** Second order one dimensional telegraph equation, Cubic trigonometric B-spline basis functions, Cubic trigonometric B-spline collocation method, Stability

## 1. Introduction
### 1.1. Problem

Consider the second order one-dimensional hyperbolic telegraph equation ("the telegraph equation"), given by

$$\frac{\partial^2 u}{\partial t^2}(x,t) + 2\alpha \frac{\partial u}{\partial t}(x,t) + \beta^2 u(x,t) = \frac{\partial^2 u}{\partial x^2}(x,t) + q(x,t) \qquad a \leq x \leq b, \quad t \geq 0 \qquad (1)$$

with initial conditions

$$u(x,0) = g_1(x), \qquad u_t(x,0) = g_2(x), \qquad a \leq x \leq b \qquad (2)$$

and the following two types of boundary conditions
1. Dirichlet boundary conditions

$$u(a,t) = f_1(t), \qquad u(b,t) = f_2(t), \qquad t \geq 0 \qquad (3)$$

2. Neumann boundary conditions

$$u_x(a,t) = w_1(t), \qquad u_x(b,t) = w_2(t), \qquad t \geq 0 \qquad (4)$$

### 1.2. Applications

The study of electric signal in a transmission line, dispersive wave propagation, pulsating blood flow in arteries and random motion of bugs along a hedge are amongst a host of physical and biological phenomena which can be described by the telegraph equation (1). Details of the above mentioned phenomena and other phenomena which can be described by the telegraph equation (1) can be found in [1-4]. Clearly the equation and its solution is of importance in many areas of applications.

*Corresponding Email address: m.abbas@uos.edu.pk

### 1.3. Literature review

Several numerical methods have been developed to solve the telegraph equation subject to Dirichlet boundary conditions and the references are in [2, 5-8]. In [9], two semi-discretization methods based on quartic splines function have been developed to solve the telegraph equations. A class of unconditionally stable finite difference schemes constructed with the help of quartic splines functions has been developed by Liu and Liu [10] for the solution of the telegraph equation. Further several numerical methods have been developed by Dehghan [11-12] in collaboration with different authors. These include the thin plate splines radial basis functions (RBF) for the numerical solution of the telegraph equation [11] and high-order compact finite difference method to solve the telegraph equation [12]. Further details on other numerical methods including interpolating scaling functions [13], radial basis functions [14], quartic B-spline collocation method (QuBSM) [15], cubic B-spline collocation method (CuBSM) [16-17] for the solution of the telegraph equation subject to Dirichlet boundary conditions are in the literature. Thus many numerical methods have been developed to solve the telegraph equation (1) with Dirichlet boundary conditions.

Some numerical methods have been developed for numerical solution of the telegraph equation with Neumann boundary conditions. These include methods by Dehghan and Ghesmati [1] who constructed a dual reciprocity boundary integral equation (DRBIE) method in which cubic-radial basis function (C-RBF), thin plate-spline- radial basis function (TPS-RBF) and linear radial basis functions (L-RBF) are utilized for the numerical solution of the telegraph equation with Neumann boundary conditions. Liu and Liu [18] have developed a compact difference unconditionally stable scheme (CDS) to solve the telegraph equation with Neumann boundary conditions. Further, Mittal and Bhatia [19] have developed a technique based on collocation of cubic B-spline collocation method (CuBSM) for solving the telegraph equation with Neumann boundary conditions.

The trigonometric B-spline collocation method has attracted attention in the literature and has been used for the numerical solutions of several linear and nonlinear partial differential equations [20-25]. The trigonometric B-splines have many geometric properties like local support, smoothness, and capability of handling local phenomena. There properties make trigonometric B-spline appropriate to solve linear and nonlinear partial differential equations easily and effortlessly. Fyfe [26] found that the spline method is better than the usual finite difference scheme because it has the flexibility to obtain the solution at any point in the domain with greater accuracy. The trigonometric B-spline produced more accurate results for linear and non-linear initial boundary value problems as compared to traditional B-spline functions [27-28].

In our work, a numerical collocation finite difference technique based on cubic trigonometric B-spline is presented for the solution of telegraph equation (1) with initial conditions in equation (2) and different two types of boundary conditions in equations (3) and (4). Several studies have been carried out as the ordinary B-spline collocation methods to solve the proposed problem subject to different types of boundary conditions but not with cubic trigonometric B-spline collocation method. A usual finite difference scheme is applied to discretize the time derivative while cubic trigonometric B-spline is utilized as an interpolating function in the space dimension. The proposed method is unconditionally stable over $0.5 \leq \theta \leq 1$ and this is proved by von Neumann approach. The feasibility of the method is shown by test problems and the approximated solutions are found to be in good agreement with the exact solutions. The proposed method is superior to C-RBF [1], TPS-RBF [1], L-RBF [1], RBF [11], QuBSM [15], CDS [18], CuBSM [16, 19] due to smaller storage and CPU time in seconds.

### 1.4. Outlines of current paper

The outline of this paper is as follows: In section 2, the cubic trigonometric B-spline collocation method is explained. In section 3, numerical solution of proposed problem (1) is discussed. In section 4, the stability of proposed method is investigated. In section 5, the results of numerical experiments are presented and compared with exact solutions and some previous methods. Finally, in section 6, the conclusion of this study is given.

## 2. Description of new Trigonometric B-spline method

In this approach, the space derivatives are approximated by using cubic trigonometric B-spline method (CuTBSM). A mesh $\Omega$ which is equally divided by knots $x_i$ into $N$ subintervals $[x_i, x_{i+1}], i = 0,1,2,...,N-1$ such that, $\Omega: a = x_0 < x_1 < ... < x_N = b$ is used. For the telegraph equation (1), an approximate solution using collocation method with cubic trigonometric B-spline is obtained in the form [27-28]

$$U(x,t) = \sum_{i=-3}^{N-1} C_i(t) TB_i(x) \tag{5}$$

where $C_i(t)$ are to be calculated for the approximated solutions $u(x,t)$ to the exact solutions $u_{exc}(x,t)$, at the point $(x_i, t_j)$. A $C^2$ piecewise cubic trigonometric B-spline basis functions $TB_i(x)$ over the uniform mesh can be defined as [20-22]

$$TB_i(x) = \frac{1}{\omega} \begin{cases} \xi^3(x_i), & x \in [x_i, x_{i+1}] \\ \xi(x_i)\left(\xi(x_i)\zeta(x_{i+2}) + \zeta(x_{i+3})\xi(x_{i+1})\right) + \zeta(x_{i+4})\xi^2(x_{i+1}), & x \in [x_{i+1}, x_{i+2}] \\ \zeta(x_{i+4})\left(\xi(x_{i+1})\zeta(x_{i+3}) + \zeta(x_{i+4})\xi(x_{i+2})\right) + \xi(x_i)\zeta^2(x_{i+3}), & x \in [x_{i+2}, x_{i+3}] \\ \zeta^3(x_{i+4}), & x \in [x_{i+3}, x_{i+4}] \end{cases} \tag{6}$$

where,

$$\xi(x_i) = \sin\left(\frac{x - x_i}{2}\right), \; \zeta(x_i) = \sin\left(\frac{x_i - x}{2}\right), \; \omega = \sin\left(\frac{h}{2}\right)\sin(h)\sin\left(\frac{3h}{2}\right)$$

where $h = (b-a)/N$.

The approximations $U_i^j$ at the point $(x_i, t_j)$ over subinterval $[x_i, x_{i+1}]$ can be defined as

$$U_i^j = \sum_{k=i-3}^{i-1} C_k^j TB_k(x) \tag{7}$$

The values of $TB_i(x)$ and its derivatives at knots are required to obtain the approximate solutions and these derivatives are recorded in Table 1.

Table. 1: Values $TB_i(x)$ and its derivatives.

|        | $x_i$ | $x_{i+1}$ | $x_{i+2}$ | $x_{i+3}$ | $x_{i+4}$ |
|--------|-------|-----------|-----------|-----------|-----------|
| $TB_i$ | 0     | $a_1$     | $a_2$     | $a_1$     | 0         |
| $TB_i'$ | 0    | $a_3$     | 0         | $a_4$     | 0         |
| $TB_i''$ | 0   | $a_5$     | $a_6$     | $a_5$     | 0         |

where

$$a_1 = \frac{\sin^2\left(\frac{h}{2}\right)}{\sin(h)\sin\left(\frac{3h}{2}\right)}, \quad a_2 = \frac{2}{1+2\cos(h)}, \quad a_3 = -\frac{3}{4\sin\left(\frac{3h}{2}\right)}, \quad a_4 = \frac{3}{4\sin\left(\frac{3h}{2}\right)}$$

$$a_5 = \frac{3(1+3\cos(h))}{16\sin^2\left(\frac{h}{2}\right)\left(2\cos\left(\frac{h}{2}\right)+\cos\left(\frac{3h}{2}\right)\right)}, \quad a_6 = -\frac{3\cos^2\left(\frac{h}{2}\right)}{2\sin^2\left(\frac{h}{2}\right)(1+2\cos(h))}$$

From (5) and (6), the values at the knots of $U_i^j$ and their derivatives up to second order are calculated in the terms of time parameters $C_i^j$ as

$$\begin{cases} U_i^j = a_1 C_{i-3}^j + a_2 C_{i-2}^j + a_1 C_{i-1}^j \\ (U_x)_i^j = a_3 C_{i-3}^j + a_4 C_{i-1}^j \\ (U_{xx})_i^j = a_5 C_{i-3}^j + a_6 C_{i-2}^j + a_5 C_{i-1}^j \end{cases} \quad (8)$$

The equation (5) and boundary conditions given in (3) and (4) are used to obtain the approximate solution at end points of the mesh as

$$\begin{cases} U(x_0, t_{j+1}) = a_1 C_{-3} + a_2 C_{-2} + a_1 C_{-1} = f_1(t_{j+1}) \\ U(x_N, t_{j+1}) = a_1 C_{N-3} + a_2 C_{N-2} + a_1 C_{N-1} = f_2(t_{j+1}) \end{cases} \quad (9)$$

and

$$\begin{cases} U_x(x_0, t_{j+1}) = a_3 C_{-3} + a_4 C_{-1} = w_1(t_{j+1}) \\ U_x(x_N, t_{j+1}) = a_3 C_{N-3} + a_4 C_{N-1} = w_2(t_{j+1}) \end{cases} \quad (10)$$

## 3. Numerical Solution of telegraph equation

In this section, a numerical solution of telegraph equation (1) is obtained using collocation approach based on cubic trigonometric basis functions. The discretization in time derivative is obtained by forward finite difference scheme and $\theta$ weighted scheme applied to problem (1) to obtain a tri-diagonal of linear equations. The proposed $\theta$ weighted scheme is closely related to the accuracy of the method and numerical stability. A uniform mesh $\Omega$ with grid points $(x_i, t_j)$ to discretize the grid region $\Delta = [a,b] \times [0,T]$ with $x_i = a + ih$, $i = 0,1,2,...,N$ and $t_j = j\Delta t$, $j = 0,1,2,3,...,M$, is used $T = M\Delta t$. The quantities $h$ and $\Delta t$ are mesh space size and time step size respectively. Using $\theta$ weighted technique, the approximations for the solutions of telegraph equation (1) at $t_{j+1}$ th time level can be given by as [21]

$$(U_{tt})_i^j + 2\alpha (U_t)_i^j = \theta g_i^{j+1} + (1-\theta) g_i^j + q(x_i, t_j) \quad \theta \in [0,1] \quad (11)$$

where $g_i^j = (U_{xx})_i^j - \beta^2 U_i^j$ and the subscripts $j$ and $j+1$ are successive time levels, $j = 0,1,2,....M$. By using the central finite difference discretization of the time derivatives and rearranging the equations (11), we obtain

$$\frac{U_i^{j+1} - 2U_i^j + U_i^{j-1}}{\Delta t^2} + 2\alpha \frac{U_i^{j+1} - U_i^j}{\Delta t} = (1-\theta) g_i^j + \theta g_i^{j+1} + q(x_i, t_j) \quad (12)$$

The equation (12) yields it as

$$(1 + 2\alpha k) U_i^{j+1} - k^2 \theta g_i^{j+1} = 2(1 + \alpha k) U_i^j + k^2 (1-\theta) g_i^j - U_i^{j-1} + k^2 q(x_i, t_j) \quad (13)$$

where $k = \Delta t$ is the time step. It is noted that the system becomes an explicit scheme when $\theta = 0$, a fully implicit scheme when $\theta = 1$, and a Crank-Nicolson scheme when $\theta = 1/2$ [20-21]. Hence, (13) becomes

$$(1+2\alpha k)U_i^{j+1} - k^2\theta\left((U_{xx})_i^{j+1} - \beta^2 U_i^{j+1}\right) = 2(1+\alpha k)U_i^j + k^2(1-\theta)\left((U_{xx})_i^j - \beta^2 U_i^j\right) - U_i^{j-1} + k^2 q(x_i, t_j)$$
(14)

The initial condition (2) is substituted into last term of equation (14) for computing $C^1$. By central difference approximation,
$$U_i^{-1} = U_i^1 - 2k\, g_2(x_i) \tag{15}$$

After that, the system thus obtained for $j \geq 1$ on simplifying (14) after using (8) consists of $N+1$ linear equations in $N+3$ unknowns $C^{j+1} = \left(C_{-3}^{j+1}, C_{-2}^{j+1}, C_{-1}^{j+1}, \ldots, C_{N-1}^{j+1}\right)$ at the time level $t = t_{j+1}$. The boundary conditions given in equations (9) or (10) are used for two additional linear equations to obtain a unique solution of the resulting system. Thus, the system becomes a matrix system of dimension $(N+3) \times (N+3)$ which is a tri-diagonal system that can be solved by the Thomas Algorithm [29-33].

## 2.1 Initial state

After the initial vectors $C^0$ have been computed from the initial conditions, the approximate solutions $U_i^{j+1}$ at a particular time level can be calculated repeatedly by solving the recurrence relation (14) [20-21]. $C^0$ can be obtained from the initial and boundary values of the derivatives of the initial condition as follows[21]:

$$\begin{cases} (U_i^0)_x = g_1'(x_i), & i = 0 \\ U_i^0 = g_1(x_i), & i = 0,1,2,\ldots,N \\ (U_i^0)_x = g_1'(x_i), & i = N \end{cases} \tag{16}$$

Thus the equations (16) yield a $(N+3) \times (N+3)$ matrix system for which the solution can be computed by the use of the Thomas algorithm.

## 3. Stability of proposed method

In this section, the von Neumann stability method is applied to investigate the stability of the proposed scheme. Such an approach has been used by many researchers [20-22, 34]. Substituting the approximate solution $U(x,t)$, their derivatives at the knots with $q(x,t) = 0$ ([35], chapter 9), into equation (14) yields a difference equation with variables $C_m$ given by

$$\begin{cases} \left((1+2\alpha k + k^2\theta\beta^2)a_1 - k^2\theta a_5\right)C_{m-3}^{j+1} + \left((1+2\alpha k + k^2\theta\beta^2)a_2 - k^2\theta a_6\right)C_{m-2}^{j+1} \\ + \left((1+2\alpha k + k^2\theta\beta^2)a_1 - k^2\theta a_5\right)C_{m-1}^{j+1} \\ = \left((2+2\alpha k - (1-\theta)k^2\beta^2)a_1 + (1-\theta)k^2 a_5\right)C_{m-3}^j + \left((2+2\alpha k - (1-\theta)k^2\beta^2)a_2 + (1-\theta)k^2 a_6\right)C_{m-2}^j \\ + \left((2+2\alpha k - (1-\theta)k^2\beta^2)a_1 + (1-\theta)k^2 a_5\right)C_{m-1}^j - \left(a_1 C_{m-3}^{j-1} + a_2 C_{m-2}^{j-1} + a_1 C_{m-1}^{j-1}\right) \end{cases}$$
(17)

Simplifying it leads to
$$w_1 C_{m-3}^{j+1} + w_2 C_{m-2}^{j+1} + w_1 C_{m-1}^{j+1} = w_3 C_{m-3}^j + w_4 C_{m-2}^j + w_3 C_{m-1}^j - a_1 C_{m-3}^{j-1} - a_2 C_{m-2}^{j-1} - a_1 C_{m-1}^{j-1} \tag{18}$$
where
$$\begin{aligned} w_1 &= (1+2\alpha k + k^2\theta\beta^2)a_1 - k^2\theta a_5, \\ w_2 &= (1+2\alpha k + k^2\theta\beta^2)a_2 - k^2\theta a_6, \\ w_3 &= (2+2\alpha k - (1-\theta)k^2\beta^2)a_1 + (1-\theta)k^2 a_5, \\ w_4 &= (2+2\alpha k - (1-\theta)k^2\beta^2)a_2 + (1-\theta)k^2 a_6 \end{aligned} \tag{19}$$

Now on inserting the trial solutions (one Fourier mode out of the full solution) at a given point $x_m$, $C_m^j = \delta^j \exp(im\eta h)$ into equation (18) and rearranging the equations, $\eta$ is the mode number, $h$ is the element size and $i^2 = -1$, we obtain

$$w_1\delta^{j+1}e^{i\eta(m-3)h} + w_2\delta^{j+1}e^{i\eta(m-2)h} + w_1\delta^{j+1}e^{i\eta(m-1)h}$$
$$= w_3\delta^{j}e^{i\eta(m-3)h} + w_4\delta^{j}e^{i\eta(m-2)h} + w_3\delta^{j}e^{i\eta(m-1)h} - a_1\delta^{j-1}e^{i\eta(m-3)h} - a_2\delta^{j-1}e^{i\eta(m-2)h} - a_1\delta^{j-1}e^{i\eta(m-1)h} \quad (20)$$

Dividing equation (20) by $\delta^{j-1}e^{i\eta(m-2)h}$ and rearranging, we obtain

$$\delta^2\left(w_2 + 2w_1\cos(\eta h)\right) - \delta\left(w_4 + 2w_3\cos(\eta h)\right) + \left(a_2 + 2a_1\cos(\eta h)\right) = 0 \quad (21)$$

The wave number is given as

$$\eta = \frac{2\pi}{\lambda} \quad (22)$$

where $\lambda$ is the wave length.
Let

$$N = \frac{\lambda}{h}$$

which represents the number of grid interval over one wavelength. Then the equation (22) can be rearranged to the form [35]

$$\varphi = \eta h = \frac{2\pi}{N} \quad (23)$$

where $\varphi = \eta h$ is dimensionless wave number. As the shortest waves represented at the considered grid points have wavelength $2h$, whereas the longest ones tend to infinity, then $2 \leq N \leq \infty$ implies that $0 \leq \varphi \leq \pi$ [35].

Let
$A = w_2 + 2w_1\cos(\varphi)$
$B = w_4 + 2w_3\cos(\varphi)$
$C = a_2 + 2a_1\cos(\varphi)$

Then the equation (21) yields

$$A\delta^2 - B\delta + C = 0 \quad (24)$$

Applying the Routh-Hurwitz criterion [34] on equation (24), the necessary and sufficient conditions for equation (14) to be unconditionally stable as follows:

Consider the transformation $\delta = \frac{1+\xi}{1-\xi}$ and simplifying the equation (24) becomes as

$$(A+B+C)\xi^2 + 2(A-C)\xi + (A-B+C) = 0 \quad (25)$$

The unconditionally stability condition $|\xi| \leq 1$ under the following necessary and sufficient conditions

$$A+B+C \geq 0,\ A-C \geq 0,\ A-B+C \geq 0 \quad (26)$$

$$A+B+C = (w_2 + w_4 + a_2) + 2(w_1 + w_3 + a_1)\cos(\varphi)$$
$$A-B+C = (w_2 - w_4 + a_2) + 2(w_1 - w_3 + a_1)\cos(\varphi) \quad (27)$$
$$A-C = (w_2 - a_2) + 2(w_1 - a_1)\cos(\varphi)$$

Since $\varphi$ ranges from 0 to $\pi$, then inequalities (26) can be verify for its extreme values only [35]. Setting $\varphi = \pi$, the values of $w_i, i = 1,2,3,4$ and $a_i, i = 1,2$, it can be easily proved that

$$A+B+C = \left(16(1+k\alpha) + 2k^2(3+4\beta^2)(-1+2\theta)\sin^2\left(\frac{h}{4}\right)\right)\cosec^2\left(\frac{h}{4}\right) \geq 0 \quad (28)$$

The inequality given in equation (28) satisfy if $-1 + 2\theta \geq 0 \Rightarrow \theta \geq \frac{1}{2}$.

$$A-B+C = k^2\cosec^2\left(\frac{h}{4}\right)\left(6\cos^2\left(\frac{h}{4}\right) + 8\beta^2\sin^2\left(\frac{h}{4}\right)\right) \geq 0 \quad (29)$$

$$A - C = k\cos ec^2\left(\frac{h}{2}\right)\left(6k\theta\cos^2\left(\frac{h}{4}\right) + 2(8\alpha + 4k\beta^2\theta)\sin^2\left(\frac{h}{4}\right)\right) \geq 0 \qquad (30)$$

Thus the proposed scheme for telegraph equation is unconditionally stable in the region $0.5 \leq \theta \leq 1$ without any restriction on grid size and time step size but $h$ should be chosen in such a way that the accuracy of the scheme is not degraded.

## 4. Numerical Experiments

This section presents some numerical results of the hyperbolic telegraph equation (1) with initial (2) and boundary conditions (3) or (4). To test the accuracy of proposed method, several numerical experiments for different values of $\alpha$ and $\beta$ are given in this section with $L_\infty$, $L_2$ and root mean square (RMS) errors are calculated by

$$L_\infty = \|u_{exc} - U_N\|_\infty = \max_j |u_j - (U_N)_j|$$

$$L_2 = \|u_{exc} - U_N\|_2 = \sqrt{h\sum_{j=o}^{N}|u_j - (U_N)_j|^2}$$

$$\text{RMS} = \sqrt{\frac{\sum_{j=o}^{N}|u_j - (U_N)_j|^2}{N+1}}.$$

We compare the numerical solutions obtained by cubic trigonometric B-spline collocation method for telegraph equation (1) with known exact solutions and those numerical methods in the literature. We carry out (14) by the proposed method and Intel ®Core ™ i5-2410M CPU@2.30 GHz with 8GB RAM and 64-bit operating system (Windows 7). The numerical implementation is carried out in Mathematica 9. Numerical results are computed by cubic trigonometric B-spline collocation method for the telegraph equation (1) at different time levels with smaller storage and CPU time which are tabulated in different Tables. All Figureures are drawn at the value of weighting parameter $\theta = 0.5$.

**Problem 1.**

Consider the following particular case of equation (1) in the domain $[0, \pi]$ with $\alpha = 2, \beta = \sqrt{2}$ [11, 16]

$$\frac{\partial^2 u}{\partial t^2}(x,t) + 4\frac{\partial u}{\partial t}(x,t) + 2u(x,t) = \frac{\partial^2 u}{\partial x^2}(x,t) + q(x,t) \qquad 0 \leq x \leq \pi, \quad t \geq 0$$

subject to the following initial and boundary conditions

$$\begin{cases} u(x,t=0) = \sin(x), & \frac{\partial u}{\partial t}(x,t=0) = -\sin(x) \\ u(x=0,t) = 0, & u(x=\pi,t) = 0 \end{cases}$$

where $q(x,t) = -2e^{-t}\sin(x)$. The exact solution of this problem is $u_{exc}(x,t) = e^{-t}\sin(x)$. The proposed method is applied to calculate the numerical solutions of the telegraph equation (1)-(3) with $h = 0.02$, $\Delta t = 0.0001$ at different time levels. The absolute errors ($L_\infty$) and relative error ($L_2$) at weighting parameter $\theta = 0.5$, different time levels and also CPU time in second, are reported in Table 2. It can be concluded that our results are more accurate as compared to results obtained by Dehghan and Shokri [11] and Mittal and Bhatia [16]. In Table 3 and Figure 1, we report the absolute errors, relative errors and RMS for $h = 0.02$, $\Delta t = 0.01$ at different time levels with CPU time with different values of weighting parameter $\theta$ due to the purpose of comparison with existing methods. The numerical results of this problem are in good agreement with exact solution and are more accurate than cubic B-spline collocation method [22]. Figure 2 depicts the graphs of comparison between exact and numerical solutions at time levels $t = 1, 2, 3$ with $h = 0.02, \Delta t = 0.01$. Figure 3 shows the space-time graph of exact and approximate solutions at $t = 3$ with $h = 0.02, \Delta t = 0.01$.

Table 2: Relative errors, maximum errors and CPU time of Problem 1 with $\Delta t = 0.0001, h = 0.02$

Table 3: Relative errors, maximum errors, RMS, values of $\theta$ and CPU time of Problem 1 with $\Delta t = 0.01, h = 0.02$

Figure 1: Error graph of Problem 1 at different time levels with $h = 0.02, \Delta t = 0.01$

Figure 2: Comparison of numerical and exact solution of Problem 1 at different time levels with $h = 0.02, \Delta t = 0.01$

Figure 3: Space-time graph of Problem 1 at $T = 3$ with $h = 0.02, \Delta t = 0.01$

**Problem 2.**

In this problem, we consider the telegraph equation (1) in the domain $[0, 2]$ with $\alpha = 10, \beta = 5$ [15-16]

$$\frac{\partial^2 u}{\partial t^2}(x,t) + 20 \frac{\partial u}{\partial t}(x,t) + 25 u(x,t) = \frac{\partial^2 u}{\partial x^2}(x,t) + q(x,t) \qquad 0 \leq x \leq 2, \quad t \geq 0$$

subject to the following initial and boundary conditions

$$\begin{cases} u(x,0) = \tan(x/2), & \frac{\partial u}{\partial t}(x,0) = \frac{1}{2}\left(1 + \tan^2(x/2)\right) \\ u(0,t) = \tan(t/2), & u(2,t) = \tan((2+t)/2) \end{cases}$$

and function $q(x,t) = \alpha\left(1 + \tan^2((x+t)/2)\right) + \beta^2 \tan((x+t)/2)$. The exact solution of this equation is $u_{exc}(x,t) = \tan((x+t)/2)$. In this problem we take $L = 2$, $h = 0.02$ and two values of time step size $k = 0.0001$ and $k = 0.001$ due to the purpose of comparison with existing methods. In Table 4 we report the absolute errors and relative errors of this problem using present method at different time levels and different values of weighting parameter $\theta$. In Table 5, we also recorded the the absolute errors and relative errors at different time levels for $h = 0.001, k = 0.001$ and concluded that our rresults are more accurate than Dosti and Nazemi [15] and Mittal and Bhatia [16]. Figure 4 illustrates the comparison of exact solution with approximate solution of this problem at various time levels and different values. In Figure 5 we show the space-time graph of approximate at time $t = 1.0$

Table 4: $L_2$, $L_\infty$ errors of Problem 2 at different time levels with $h = 0.02$

Table 5: Relative errors and maximum errors of Problem 2 with $\Delta t = h = 0.001$

Figure 4: Numerical and exact solutions of Problem 2 at different time levels with $h = 0.02, \Delta t = 0.001$

Figure 5: Space-time graph of approximate solution of Problem 2 at $t=1.0$ with $h=0.02$

**Problem 3**

We consider the telegraph equation (1) in the domain $[0,1]$ with $\alpha=0.5, \beta=1.0$ [11, 16]

$$\frac{\partial^2 u}{\partial t^2}(x,t)+\frac{\partial u}{\partial t}(x,t)+u(x,t)=\frac{\partial^2 u}{\partial x^2}(x,t)+q(x,t)$$

subject to the following initial and boundary conditions

$$\left\{u(x,0)=0,\ \frac{\partial u}{\partial t}(x,0)=0\quad 0\le x\le 1\ \&\ u(0,t)=0,\ u(1,t)=0\quad t\ge 0\right.$$

and $q(x,t)=\left(2-2t+t^2\right)\left(x-x^2\right)e^{-t}+2t^2 e^{-t}$. The exact solution of this problem is $u_{exc}(x,t)=\left(x-x^2\right)t^2 e^{-t}$. The absolute errors, relative errors and CPU time in seconds is shown in Table 6 with $\Delta t=0.001,\ h=0.01$. Numerical results are compared with the obtained results in Dehghan and Shokri [11] and Mittal and Bhatia [16]. It can be concluded that the numerical solutions obtained by our method are good in comparison with Dehghan and Shokri [11] and Mittal and Bhatia [16]. The graph of exact and numerical solutions at $t=1,2,3,4,5$ is shown in Figure 6 and the space time graph of numerical solution up to t = 5 is presented in Figure 7.

Table 6: Relative errors, maximum errors and CPU time of Problem 3 with $\Delta t=0.001,\ h=0.01$

Figure 6: Comparison of Numerical and exact solutions of Problem 3 at different time levels with $h=0.01, \Delta t=0.001$

Figure 7: Space-time graph of approximate solution of Problem 3 at $t=5.0$ with $h=0.01$

**Example 4**

Consider the telegraph equation (1) in the domain $[0,1]$ and $\alpha=6, \beta=2$ [15-16]

$$\frac{\partial^2 u}{\partial t^2}(x,t)+12\frac{\partial u}{\partial t}(x,t)+4u(x,t)=\frac{\partial^2 u}{\partial x^2}(x,t)+q(x,t)$$

with following initial and boundary conditions

$$\begin{cases} u(x,0)=\sin(x), & \frac{\partial u}{\partial t}(x,0)=0 & 0\le x\le 1 \\ u(0,t)=0, & u(1,t)=\cos(t)\sin(1) & t\ge 0 \end{cases}$$

and $q(x,t)=-2\alpha\sin(t)\sin(x)+\beta^2\cos(t)\sin(x)$.

The exact solution of this problem is $u(x,t)=\cos(t)\sin(x)$.

The efficiency can be noted from Table 7 by using $L_2, L_\infty$ and RMS errors with $\Delta t=0.0001, h=0.01$. In Table 8 we also reported the absolute errors and relative errors at different time levels for different values $h$ with $\Delta t=0.001$ and the numerical results are compared with those of Dosti and Nazemi [15] and Mittal and Bhatia [16]. We found that our numerical results are comparable to that of QuBSM [15] and CuBSM [16] in terms of $L_2, L_\infty$ errors. Figure. 8 presents the comparison of numerical and exact solutions for different time levels with $\Delta t=0.001, h=0.01$. The space–time graph of numerical solution at $t=1.0$ is presented in Figure. 9.

Table 7: Relative errors, maximum errors and RMS errors of Problem 4 with $\Delta t=0.001, h=0.01$

Table 8: $L_2$, $L_\infty$ and RMS errors of Problem 4 at different time levels with $k = 0.001$

Figure 8: Comparison of Numerical and exact solutions of Problem4 at different time levels with $h = 0.01, \Delta t = 0.001$

Figure 9: Space-time graph of approximate solution of Problem 4 at $t = 5.0$ with $h = 0.01$

**Example 5**
Consider the following particular case of second order one dimensional equation (1) over the region $[0, 2\pi] \times [0, 3]$ with $\alpha = 4, \beta = 2$ [1, 18-19]

$$\frac{\partial^2 u}{\partial t^2}(x,t) + 8\frac{\partial u}{\partial t}(x,t) + 4u(x,t) = \frac{\partial^2 u}{\partial x^2}(x,t) + q(x,t) \qquad 0 \leq x \leq 2\pi, \quad t \geq 0$$

subject to the following initial and Neumann boundary conditions

$$\begin{cases} u(x, t=0) = \sin(x), & \frac{\partial u}{\partial t}(x, t=0) = -\sin(x) \\ u_x(0,t) = e^{-t}, & u_x(2\pi, t) = e^{-t} \end{cases}$$

where $q(x,t) = -2e^{-t}\sin(x)$. The exact solution of this problem is $u_{exc}(x,t) = e^{-t}\sin(x)$. The proposed method is applied to calculate the numerical solutions of telegraph equation (1)-(2) and (4) at $t = 1, 2, 3$ with $\Delta t = 0.01$ and different values of $h$. The absolute errors, relative errors and RMS errors at different values of weighting parameter $\theta$ and also CPU time in second, are reported in Table 9. It can be concluded that our results are more accurate as compared to results obtained by three radial basis functions schemes such as Cubic-RBF (CRBF) [1], Thin Plate Spline RBF (TPS-RBF) [1], Linear RBF (L-RBF) [1], CDS [18] and CuBSM [19]. Figure 10 depicts the errors of proposed method at different values of $h$. The numerical results of this problem are also in good agreement with exact solution. Figure 11 shows the space-time graph of approximate solutions at $t = 3$ with $h = 0.05, \Delta t = 0.01$.

Table 9: $L_2$, $L_\infty$, RMS errors and CPU (s) of Problem 5 at different time levels with $k = 0.01$

Figure 10: Error graph of Problem 5 at different values of $h$ at $t = 3.0$ with $\Delta t = 0.01$

Figure 11: Space-time graph of the approximate solution of Problem 5 up to $t = 3.0$ with $k = 0.01$ and $h = 0.05$

**6. Conclusions**
This paper has investigated the application of cubic trigonometric B-spline collocation method to find the numerical solution of the telegraph equation with initial condition and Dirichlet as well as Neumann's type boundary conditions. A usual finite difference approach is used to discretize the time derivatives. The cubic trigonometric B-spline is used for interpolating the solutions at each time. The numerical results shown in Tables 2-9 and Figures1-11 indicate the reliability of results obtained. The obtained solution to the telegraph equation for various time levels has been compared with the exact solution and existing methods by calculating $L_\infty$, $L_2$ and RMS errors. The comparison indicated improved accuracy compared to C-RBF [1], TPS-RBF [1], L-RBF [1], RBF [11], QuBSM [15], CDS [18], CuBSM [16, 19].